\documentclass[12pt]{article}

\title{ Glimpses of the  Octonions  and Quaternions History and Today's Applications in Quantum Physics (*)}
\author {A.Krzysztof Kwa\'sniewski \\
the Dissident \\
- relegated by Bia\l ystok University  authorities  \\
from the Institute of Computer Sciences\\
organized by  the author \\
to Faculty of Physics\\ 
ul. Lipowa 41,  15 424  Bia\l ystok, Poland\\
e-mail: kwandr@gmail.com, http://ii.uwb.edu.pl/akk/publ1.htm}

\usepackage{graphicx, graphics} 
\usepackage{amsfonts,amssymb}
\usepackage{amscd} 
\chardef\bslash=`\\ 
\begin{document}
\maketitle

\vspace{2mm}
\begin{abstract}
Before we dive into the accessibility stream of nowadays indicatory applications of octonions to computer and other sciences and to quantum physics let us focus for a while on the crucially relevant  events for today's revival on interest to nonassociativity.  Our reflections keep wandering back to  the Brahmagupta-Fibonacci Two-Square Identity and then  via the Euler Four-Square Identity up to the Degen-Graves-Cayley Eight-Square Identity.
These glimpses of history incline and invite us  to re-tell  the story  on how about one month after quaternions have been carved on the Brougham bridge octonions were discovered by John Thomas Graves (1806-1870), jurist and mathematician - a friend of William Rowan Hamilton (1805-1865).
As for today we just mention en passant quaternionic and octonionic quantum mechanics, generalization of Cauchy-Riemann equations for octonions and   Triality Principle and $G_2$ group in spinor language in a descriptive way  in order not  to daunt  non-specialists.   Relation to finite geometries is recalled and the links to the  7Stones of seven sphere , seven "imaginary"  octonions'  units  in out of the Plato's Cave Reality applications are appointed.This way  we are welcome back to primary ideas  of Heisenberg, Wheeler and other distinguished founders of quantum mechanics and quantum gravity foundations.
\end{abstract}

\vspace{0.4cm}


\vspace{0.1cm}

\noindent AMS Classification Numbers: 17A75, 51E15, 81Q60

\vspace{1cm}

\noindent (*) Lectures prepared for  International Seminar on History of Mathematics in memory of Subhash Handa on December 17-18,
2007 at Ramjas College, University of Delhi and International Conference on Advances in Mathematics:
Historical Developments and Engineering Applications from 19th to 22nd December, 2007, organized by  
Indian Society for History of Mathematics and G. B. Pant University of Agriculture and Technology, Pant Nagar,INDIA.\\

\vspace{2mm}

\section{ Finite geometries and infinite poem}
\vspace{1mm}
\noindent With this essay we add to Martin Huxley's \textit{An infinite Poem} [1]  another  strophe:

\begin{quote}

\noindent CAYLEY and  HAMILTON and John Thomas GRAVES\\
\noindent those were and are The BRAVES\\
\noindent since CAYLAY - DICKSON this is the case\\
\noindent that two quaternions make  octaven  tool\\ 
\noindent then to be used with $G_2$ group  too ...\\
\end{quote}

This note being on infinite subject has also its references grouped into two parts.
First part starts from [1] and runs up to [40] - \cite{1,2,3,4,5,6,7,8,9,10,11,12,13,14,15,16,17,18,19,20,21,22,23,24,25,26,27,28,29,30,31,32,33,34,35,36,37,38,39,40} - which might mark  this note temporal end as these are sufficiently abundant and indicative.
But life goes on and then follows an endless sequence, starting from [41].
Moreover. The number of references tells you how old I am. Temporally.

\vspace{2mm}

\noindent \textbf{I. Galois numbers } 

\noindent We use the standard convention.  Out of  the vector space $V(n+1,q)$  we derive the geometric structure $PG(n,q) $ and this will be called the projective geometry of dimension $n$  over $GF(q) $  denoting Galois field. The dimension  is to meant so that  lines have 1 dimension, planes have 2 dimensions, etc.  
\noindent Let  now $S$ be the family of all finite dimensional subspaces $V(n,q)$ of an infinite dimensional linear space over finite Galois field $GF(q)$.  Consider   these finite geometries  posets  $\langle S, \subseteq \rangle$  and let  define the  Galois number  $ G_{n,q} $  to be

\begin{center}
	$G_{n,q} = \sum_{k}^{n} {n \choose k}_q$ = \textbf{Galois number} = the number of all subspaces of $V(n,q)$
\end{center}

\noindent Then for $\langle S, \subseteq \rangle $ we awake 

$$
	G_{n,q} = \left| \left[ x,y \right] \right|
$$

\noindent for such  $x,y \in S$  that maximal chain joining $x$ with $y$  is of length  $n$.

\noindent \textbf{Galois numbers} play the crucial role in the finite geometries  poset $\langle S, \subseteq \rangle$   description  and applications including  such identities as the following one "`\textbf{q-exp  squared identity}"'

$$
	[\mathit{exp}_q(x)]^2 \equiv 
	\left( \sum_{n\geq 0} \frac{x^n}{(n_q)!} \right)^2 \equiv
	\sum_{n\geq 0} \frac{G_{n,q}x^n}{(1+q)(1+q+q^2)...(1+q+...+q^{n-1})} 
$$

\noindent \textbf{II. Lattices of subspaces and projective spaces }
\vspace{2mm}	 

\noindent In accordance with ${3 \choose 1}_2 = {3 \choose 2}_2=7 $   and      $G_{3,2} = 1+7+7+1=16$   we
characterize the adequate finite geometry  for $ n=3, q=2$   by  the Hasse diagram of $L(3.2)$  in   Fig.1.

\begin{figure}[ht]
\begin{center}
\includegraphics[width=60mm]{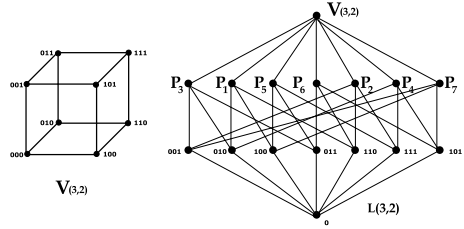}
\caption{Display of V(3.2) and L(3.2)} 
\end{center}
\end {figure} 
\vspace{2mm}

\noindent The corresponding   projective space $P(2.2)$  is then  the plane consisting of $7$  points [ see: non-zero elements - hence designating lines in $V(3.2)$ ]  and  $7$  projective lines [ see: $P_1$,. . . ,$P_7$ planes - subspaces of  $V(3.2)$ ].  It is depicted by Fig.2.

\vspace{2mm}

\begin{figure}[ht]
\begin{center}
\includegraphics[width=60mm]{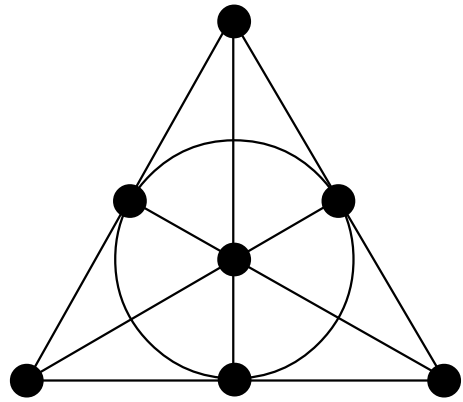}
\caption{Display of Fano Duality Plane} 
\end{center}
\end{figure}
 
\vspace{2mm}

\noindent Fig.2.  delineates the smallest  possible projective plane $P(2.2)$  known as  the  Fano plane. In accordance with the Fig.1. it contains  seven points and seven lines. These seven points are shown as  dots of intersections, while the seven projective lines are just  six line segments and a circle. Note that  one could equivalently consider the dots to be the lines while  the line segments and circle to be the points. This property is an example of the known duality property of projective planes. This means that  if the points and lines are interchanged, the resulting  set of objects  is again a projective plane. We are dealing with octonions then. See more what follows.

\noindent In accordance with ${2 \choose 1}_2 = 3 $   and      $G_{2,2} = 1+ 3+ 1= 5$   we
characterize the adequate finite geometry  for $ n= 2, q=2$   by  the Hasse diagram of $L(2.2)$  in   Fig.3.
where also the  [circle]  projective line  $P(1.2)$   is displayed.

\vspace{2mm}

\begin{figure}[ht]
\begin{center}
\includegraphics[width=60mm]{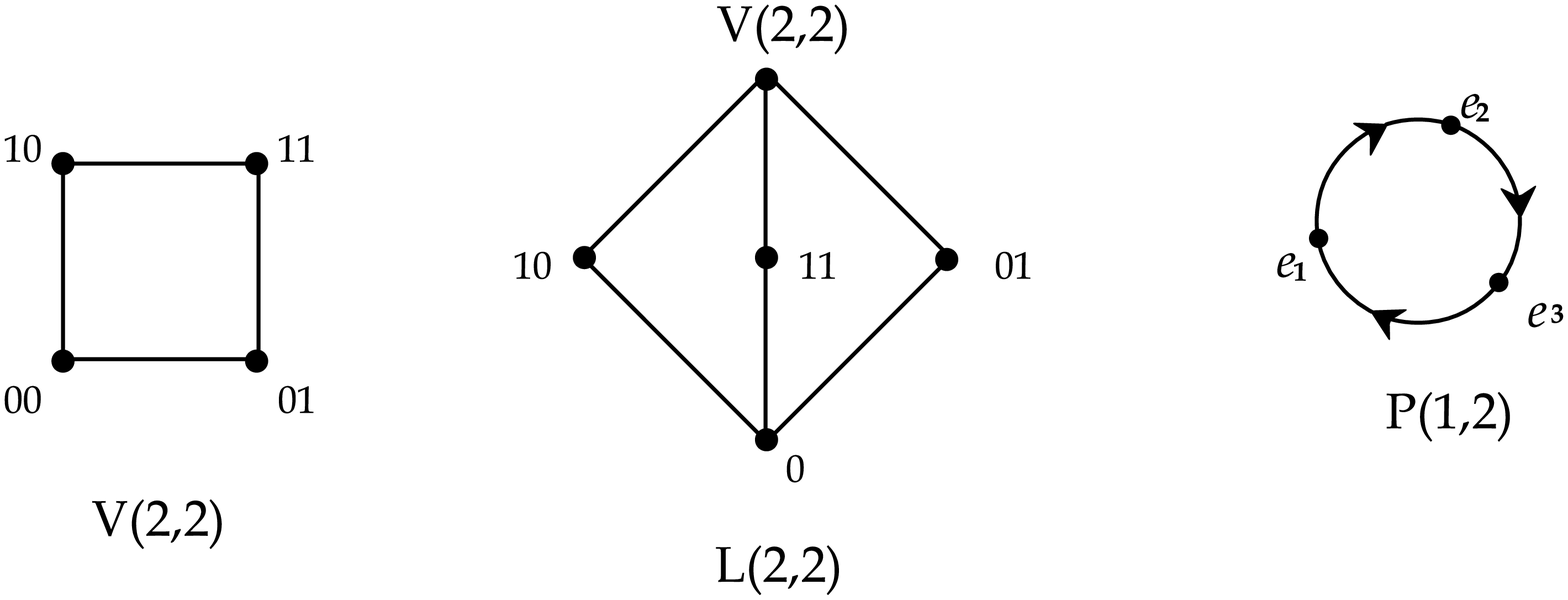}
\caption{Display of V(2.2), L(2.2) and P(1.2)} 
\end{center}
\end{figure}

\vspace{2mm}
\noindent We are dealing with quaternions then. See more what follows.
\noindent There are many ways to correlate  octonions [quaternions]  multiplication tables rules  with underlying projective 
geometry  picture [9,5,10,15]. 
Here come quite obvious examples. Using the following adjusted to Fig.1  identification 

\noindent \textbf{Projective lines }$L_i$   in $P(2.2)$,  $i = 1,2,3,4,5,6,7$, $V(3.2)$ - corresponding  Subspaces-\textbf{Planes } $P_i$\\
$$L_1   =  \left\{100,  110,   010 \right\}      \succ     P_1   = {000, 100,  110,   010}$$        
$$  L_2   =    \left\{001,  111,   110 \right\}        \succ         P_2   = \left\{000, 001,  111,   110\right\} $$
$$  L_3   =    \left\{010,  011,   001 \right\}     \succ            P_3   = \left\{000, 010,  011,   001\right\}  $$        
$$  L_4   =    \left\{010,  111,   101 \right\}        \succ         P_4   = \left\{000, 010,  111,  101\right\}     $$     
$$  L_5   =    \left\{100,  111,   011 \right\}         \succ        P_5   =\left\{000, 100,  111,  011\right\}       $$  
$$  L_6    =    \left\{011,  101,   110 \right\}        \succ        P_6   = \left\{000, 011,  101,   110\right\}        $$
$$  L_7   =     \left\{100,  101,   001 \right\}        \succ        P_7   = \left\{000, 100,  101,   001\right\}     $$

\vspace{2mm} and the following identification of octonion imaginary units
$$e_1 = 010,     e_2=100,   e_3 =110, e_4=001,  e_5= 011,  e_6 =101,  e_7= 111$$
\noindent we have

\begin{figure}[ht]
\begin{center}
\includegraphics[width=60mm]{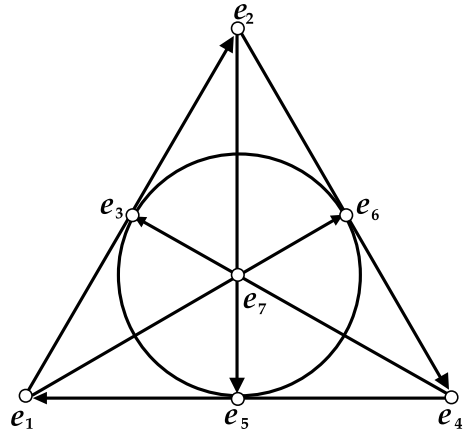}
\caption{ Fig.4  Display of octonion multiplication rules Fano Plane} 
\end{center}
\end{figure}

\vspace{2mm} 

\noindent \textbf{The rules:}
\vspace{2mm} 
$$  e_1 e_3 = e_2,  e_2 e_6 = e_4,  e_4 e_5 = e_1,  e_3 e_6 = e_5,  e_1 e_7 = e_6,
 e_2 e_7 = e_5,  e_4 e_7 = e_3$$   

\noindent Accordingly

$$ L_1 = \left\{e_1, e_2, e_3\right\}, L_2 = \left\{e_4, e_7,  e_3 \right\}, L_3 = \left\{e_1, e_5, e_4 \right\}, L_4 = \left\{e_1, e_7, e_6 \right\},$$                           
$$ L_5 = \left\{e_2, e_7, e_5 \right\}, L_6 = \left\{e_5,  e_6,  e_3 \right\}, L_7  = \left\{e_2, e_6, e_4 \right\}.$$

\noindent Using another  adjusted to Fig.1  identification  of octonion imaginary units
$$e_1= p_4=010, e_2 = p_1=100, e_3 = p_2=110, e_4 = p_3=001, $$
$$e_7= p_5 = 111,  e_5 = p_6=011,  e_6 =p_7=101 .$$

\noindent we we get the Fano Plane imported from  Joseph Malkevitch's "Finite Geometries" [12]; see Fig.5.

\begin{figure}[ht]
\begin{center}
\includegraphics[width=60mm]{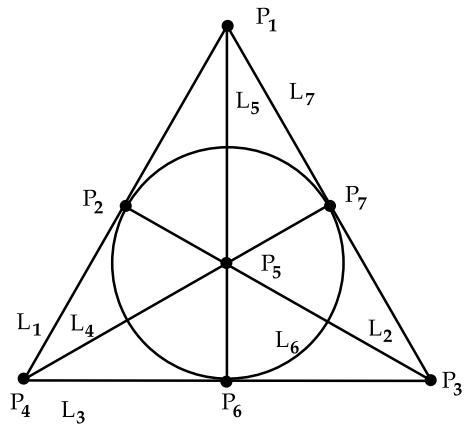}
\caption{Display of another octonion Fano Plane coding} 
\end{center}
\end{figure}

\vspace{2mm}

\noindent The still another way is the self-explanatory way of picture rule for octonions multiplication 
borrowed from  Tevian Dray's garner [TD], where the corresponding identifications with preceding 
presentations are evident. [TD]=Tevian Dray Octonions \\
http://www.physics.orst.edu/~tevian/octonions/\\, see more:
http://www.math.oregonstate.edu/~tevian/

\begin{figure}[ht]
\begin{center}
\includegraphics[width=60mm]{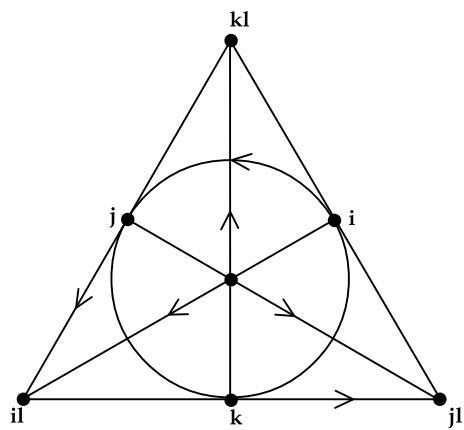}
\caption{Still another octonion rules coding} 
\end{center}
\end{figure}

\vspace{2mm}

\vspace{2mm}
\noindent Ending this finite geometry glimpses subsection we refer the reader to Online Notes in Finite Geometry http://cage.ugent.be/~aoffer/fgw/notes.html by many distinguished authors and there see Notes on Finite Geometry  
by Steven H. Cullinane  http://log24.com/notes/, where the summary Steven H. Cullinane 's work in the area of finite geometry during the years 1975 through 2006  is to be found.

\vspace{2mm}

\noindent \textbf{III. Octonions, quaternions - wide open access }
Here now follows  afar not complete selective subjectively personal overview description of the nowadays sources
\textbf{abundant with references} and links to the topics on division algebras in various branches of mathematics and other sciences 
- intriguing history included such as [2,3,4] and much more deep into the history linking us to the noted Indian mathematician 
Brahmagupta (598-c.- 665) and Bhaskara Bhaskaracharya II  (1114-1185) the lineal successor of  Brahmagupta. Bhaskara Bhaskaracharya II 
also called  Bhaskaracarya, or Bhaskara The Learned was  the leading mathematician of the 12th century on our planet. It was Him  who wrote the first  tractatus  with complete and systematic use of the decimal number system. I cannot resist- recalling my first experiences with mathemagics had been reading   the Lilavati ("The Beautiful") in Polish by Szczepan Jele\'nski. The title commemorates and refers to Bhaskara treatize.  \textit{Kalejdoskop matematyczny} by W{\l}adys{\l}aw Hugo Dionizy Steinhaus - a reward-gift from  him to me in my childhood and then Lilavati those were my fist books on the Lilavati  Mathemagics; [ Mathemagics?  see:  http://ii.uwb.edu.pl/akk/   ].

\vspace{2mm}

\noindent \textbf{Finite geometries. } Apart from the reference to finite geometries the book of Professor James Hirschfeld 
on Projective Geometries over Finite Fields is recommended for the advanced acute and persistent readers. 
Another source  is nice and useful Timothy Peil's  Survey of Geometry:\\ 
http://www.mnstate.edu/peil/geometry/IndexF/indexold.htm.

\vspace{2mm}
\noindent \textbf{ Division algebras: octonions, etc.} Here one is welcomed to the cosmic realm of John C. Baez, [http://math.ucr.edu/home/baez/README.html]  with his unsurpassed "This Week's Finds In Mathematical Physics" and the famous article [5] \textit{ The Octonions  } from which we quote the Abstract  in order to indicate genuine relevance and  scope of his ideas and "`weeks"' writings.
\noindent "The octonions are the largest of the four normed division algebras. While somewhat neglected due to their nonassociativity, they stand at the crossroads of many interesting fields of mathematics. Here we describe them and their relation to Clifford algebras and spinors, Bott periodicity, projective and Lorentzian geometry, Jordan algebras, and the exceptional Lie groups. We also touch upon their applications in quantum logic, special relativity and supersymmetry."

\noindent Also the book On Quaternions and Octonions: Their Geometry, Arithmetic, and Symmetry by John H. Conway and Derek A. Smith [6] is recommended as well as the  Geoffrey M. Dixon, book on Division Algebras : Octonions, Quaternions, Complex Numbers and the Algebraic Design of Physics [9].  A very important for this essay goals are the articles of of Titus Piezas III  and especially this on "The Degen-Graves-Cayley Eight-Square Identity"  from The Ramanujan Pages [11].  If you are still in need for short [5 pages! short] Bibliography on quaternios starting with 19-th century - then see John H. Mathews, Bibliography for Quaternions using for example password http://en.wikipedia.org/wiki/Quaternion. It starts from    \textit{Quaternions} by 
Christine Ladd  in The Analyst, Vol. 4, No. 6. (Nov., \textbf{1877}), pp. 172-174,  and ends with    
\textit{Passive Attitude Control of Flexible Spacecraft from Quaternion Measurements} by Di Gennaro S. in Journal of Optimization Theory and Applications, January \textbf{2003}, vol. \textbf{116}, no. 1, pp. 41-60.

\vspace{1mm}

\noindent For  Hamilton's Discovery of Quaternions   see excellent  paper [19]  of the Master B. L. van der Waerden.  

\vspace{1mm}

\noindent Bibliography for octonions? And quaternions? And Clifford algebras?  See 94 references in [5] including Jacques Tits papers, Pascual Jordan fundamental papers, Anthony Sudbery papers , V. S. \textbf{Varadarajan} book, Ian R. Porteous' \textit{book} [he showed \textit{it} fresh at the Catenbury 1985 conference devoted to my beloved Clifford algebras], then Feza Gürsey and Chia-Hsiung Tze papers [ I met? one of them in Paradise  [Cocoyoc in 1980]], and first of all - the famous paper by Pascual Jordan, John von Neumann, Eugene Wigner,(\textbf{1934}) [ I had met no one of them exept Wigner] however I had met all of them reading - while studying in   Wroc{\l}aw - with passion \textit{On an algebraic generalization of the quantum mechanical formalism}, (\textbf{1934}). Are there any important others' contributions ? see also [9,10,15]. And ...there is no end until The End.

\vspace{2mm}

\noindent It is quite well known how to  use the octonions  for building all five of the exceptional simple Lie algebras  including  $g_2$   intrinsic relevant to parallelizable spheres and triality principle and also interesting for quite abstract immortal theoretical physics [16,17,34].  Here the "magic" use of octonions is crucial if we are not to be overloaded by the much structured respectable and important language of all  sophisticated  Mathematics for those The Chosen. In  \textit{The Octonions Yet Another Essay } by David A. Richter one may find      among others the exceptional Lie algebra $g_2$  treatment  with overall ingenious simple use of octonions   [http://homepages.wmich.edu/~drichter/octonions.htm]  as well a s the Hamming code and Gosset's lattice  and  clear realization in octonion geometry picture of an  action by the dihedral group $D_6$  which is  the Weyl group for the exceptional 14-dimensional simple Lie algebra $g_2$.

As for the quaternion valued functions' analyticity we postpone here any attempt to master the subject of references. Let us however remark 
on the occasion that that the Fueter analyticity is a special case of Clifford analyticity. This is too big theme for glimpses. For octave Cauchy-Riemann equations we refer here to [13,14]  as we have no contact with eventual recent developments. Though we shall come back for a while to both quaternion and octonion valued functions (compare with [26,27]) when telling a bit on applications of both in  Quantum Mechanical quite recent models or projects in the almost last section with  7Stones aboard  { hence [7]).

\section{John Thomas Graves discovers Cayley numbers before Caylay}

\noindent Bhaskara in his \textit{Bijaganita} ("Seed Counting"), compiled problems from Brahmagupta and others. He filled many of the gaps in Brahmagupta's work, especially in obtaining a general solution to the \textbf{Pell} equation" . See\\ http://www.math.sfu.ca/histmath/India/12thCenturyAD/Bhaskara.html\\
also for his poetry in mathematical writings. John Pell: 1611 - 1685. Brahmagupta: 598 - 670. It is therefore around 1000 years before Pell when Brahmagupta  studied this equation with his  Brahmagupta's lemma proved by Brahmagupta in 628 AD by the method called samasa by  Indian Mathematicians.
see more in Article by: 
J J O'Connor and E F Robertson ; February 2002, MacTutor History of Mathematics, http://www-history.mcs.st-andrews.ac.uk/HistTopics/Pell.html  and references therein.

\vspace{1mm}

\noindent The similar story is about history of what we now call also after Titus Piezas III [11]  The Brahmagupta-Fibonacci Two-Square Identity  ("'just"' complex numbers) :

$$  (a^2 + b^2) (c^2 + d^2) = (ac - bd)^2 + (bc + ad)^2 . $$

\noindent At Brahmagupta and then at Fibonacci times this identity was a mystery and enigma - valid for all numbers of those times identity ?!
Well- today it is "just" norm composition under product property and Hurwitz theorem but at Brahmagupta Fibonacci times...
It is a mystery of our times. The norm of the product is equal to the product of the norms. Why only four? We now know why - anyhow it is a mystery.
The Mystery.
\noindent The Euler Four-Square Identity and the Degen-Graves-Cayley Eight-Square Identity are of the same origin as explained above: "just" norm composition under multiplication property and the Hurwitz's theorem for composition algebras.
The \textbf{Euler} Four-Square Identity \textbf{?}  Leonhard Euler (1707 - 1783) was born in Switzerland, studied under Johann Bernoulli or Ivan Bernoulli in Petersburg [20]).   Sir William Rowan Hamilton (1805 – 1865) had discovered quaternions in 1843 [19]. Hence again - Four-Square Identity was known before quaternions discovery was carved on Broome Bridge in Dublin on October 16 (1843 ) while William Rowan was being out walking along the Royal Canal  with his wife. That is mathematicians' wifes destination not only in Ireland.

\begin{figure}[ht]
\begin{center}
\includegraphics[width=75mm]{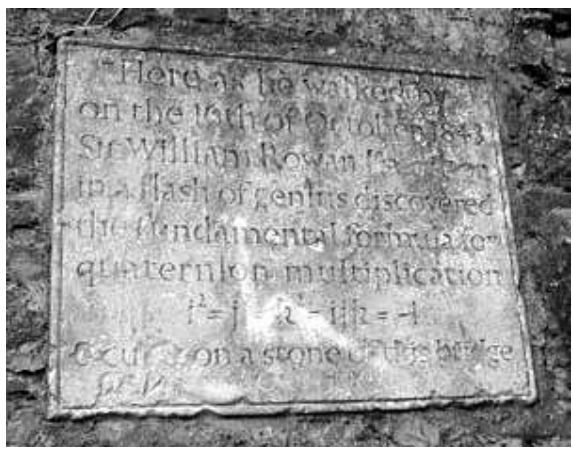}
\caption{Broome Bridge}
\end{center}
\end{figure}

\vspace{2mm}

\noindent The \textbf{Degen}-Graves-Cayley Eight-Square Identity? First discovered by Ferdinand Degen (1766-1825) the Danish mathematician around 1818. Hence again - Eiht-Square Identity was discovered before octaves.  

\noindent The Degen- \textbf{Graves}-Cayley Eight-Square Identity? It was subsequently independently rediscovered twice: in 1843 by the jurist and mathematician \textbf{John Thomas Graves} (1806-1870) and in 1845 by Arthur Cayley (1821-1895). Naurally again  the identity follows from the fact that the norm of the product of two octonions is the product of their norms. Octonions were frequently named  Cayley numbers - even in [14].
\noindent Octaves? Freudenthal's Oktaven? [21] ([21] is reference of [13,14]. I have got a copy of [21] in Russian from the late Professor Ogievietski  at Dubna Institute for Nuclear Research years, years ago.... It was rather top secret! 

\vspace{2mm}

\noindent Do you still remember? - It was  Freudenthal who had found out that the five  simple Lie groups which are out of the four infinite classical families ( exceptional groups and algebras) are related to the isometries of octaval planes [35-37]. For more see [39]; there in  P. Kainen's octonion model for physics one encounters many fascinating ideas and vague references. Read and listen (from [39] and ad rem):
\vspace{2mm}
\noindent "`In mathematics too, one finds higher-dimensional \textbf{objects casting shadows in the lower dimensions}"'. Plato's cave? "`For example, the Penrose non-periodic tiling of the plane is a projection of something in dimension at least $4$ (Katz) and the Hardy-Ramanujan formula (Andrews) shows that the number of partitions of a positive integer $n$  may be expressed in terms of the $24 k$-th roots of unity.  While integers are $0$-dimensional, the partition formula suggests that a $24$-fold symmetry is involved - as would be the case if it involves objects in $4$-dimensions."'

\vspace{2mm}

\noindent \textbf{Objects casting shadows in the lower dimensions?}

This Plato's cave idea of shadows is everywhere present in an \textbf{vedic-minded } John Archibald Wheeler's writings . Note his famous device motto 
"`\textbf{It from Bit}"'.
\noindent Wheeler wrote: "I like to think that someone will trace how the deepest thinking of \textbf{India} made its way to \textbf{Greece }and from there to the philosophy of our times." 

 "`It is curious that people like Schroedinger, Niels Bohr, Oppenheimer and \textbf{John Wheeler }are \textbf{Upanishad} scholars."' This was another glimpse from the  article - Indian Conquests of the Mind - By \textbf{Saibal Gupta.} The Statesman.org). The next glimpse is:
\noindent "One has the feeling that the thinkers of the East knew it all, and if we could only translate their answers into our language we would have the answers to all our questions." 
(source: \textit{Uncommon Wisdom} - By Fritjof Capra p. 40).

\noindent  Werner Heisenberg should be perhaps also counted an Upanishad scholar. For sure he was neoplatonic in thinking and feeling Quantum Mechanics message. Werner Karl Heisenberg (1901-1976)  German theoretical physicist was one of the leading scientists of the 20th century. Heisenberg spent some time in India as \textbf{Rabindranath Tagore}'s guest in 1929. There he got acquainted with Indian philosophy which brought him great comfort for its similarity to modern physics. 
Anyhow is it not like that with all of them:   Wheeler, J. A. and Zurek, W. H. and  Schroedinger, Niels Bohr, Oppenheimer and you and me - all we are 
\textbf{It from Bit}?  This was a kind of passionate idea of Feynman at the dusk of his life.\\
See  \textit{Simulating physics with computers }by Feynman, IJTP 21 1982 [contact !  http://www.stardrive.org/Jack/Feynman.html ], See Tony Hey  memorial  http://www.quniverse.sk/buzek/zaujimave/p257\_s.pdf.  ... \textbf{It from Bit}? 

\vspace{2mm}

\noindent Ad  \textbf{It form Bit}and my private Quantum  Plato's cave use Google password    \textbf{John Archibald Wheeler to Tachion}  and then once you are on  Home Page AKK . There in \textbf{publikacje} http://ii.uwb.edu.pl/akk/publ1.htm   see and/or download  the reference [38] here and ref. [38] there. If you click the link \textbf{John Archibald Wheeler to Tachion} right at the beginnig of my home page then you will see something that you have never seen before , something written by John Archibald Wheeler to Tachion  i.e. to me. As now I am the local dissident so I have got a resistance nick-name : Tachion  i.e. a person who causes authorities' contraventions in their past  reacting on these from their future.

\noindent Let us however continue our main story. Our story  without authorities.

\noindent Octaves?  This is a name given to octonions by Hamilton. This  story and the story on how octonions were discovered by John Thomas Graves is quite well known. Perhaps the most complete available knowledge is  John Baez's week152 [8] and links there and [19] masterpiece. Therefore instead of repeating the unsurpassed in naturality  and reliability John Baez's story  we  add something perhaps not that well known to the audience curious in the circumstances we are now involved with interest in.
Let me indicate the presence of John Thomas Graves' \textbf{brothers} in his and Sir William Rowan Hamilton's quternion-octaven life story.
Let me tell you something about \textbf{Charles Graves }  Born: 12 November (or 6 December) 1812 and Died: 1899.   
\noindent Sir William Rowan Hamilton died in 1865.  The Memorial Address was delivered by Charles Graves, who was President of the Royal Irish Academy at the time of Hamilton's death\\ 
http://www.maths.tcd.ie/pub/HistMath/People/Hamilton/Eloge/Eloge.html\\

\noindent His Éloge was delivered at the Stated Meeting of the Royal Irish Academy on the 30th November, 1865.  The Graves brothers had known Hamilton from his college days: the elder brother John Graves had been a classmate of Hamilton. Charles Graves had been a friend and a colleague of Hamilton, as Professor of Mathematics at Trinity
1829-1835 Trinity College in Dublin and at the same time  was a brother of Robert Perceval Graves who was Hamilton's first biographer. 

\vspace{2mm}

\noindent Charles Graves was elevated up to the position of Bishop of Limerick. He was a noted linguist and antiquarian, and was the grandfather of the poet Robert Graves. This is known quite well. The less known fact [I have learned about it from  Professor O.V. Viskov from Steklov Institute] - the less known fact is that the Bishop of Limerick - Charles Graves in [22] applied what is now sometimes called the Blissard umbral calculus and now also referred to Rota the operator algebra. This algebra had been practiced already in 19-th century, then the algebra   rediscovered in quantum physics as Heisenberg-Weyl algebra in 20-th century [23,24].  This is  crucial algebra of Quantum Theories including lasers. The present author in a bunch of research papers had introduced the name \textbf{Graves Heisenberg Weyl Algebra} - GHW Algebra in short, for this nowadays exploited algebra in plenty of applications. The name GHW algebra is already known to and has been e-mail discussed with the distinguished Professor  Patric D.F. Ion!

\noindent For GHW very first references and sometimes the papers to be loaded it is enough to write into Google the password: \textit{Graves Heisenberg Weyl Algebra}. One is referred also to my another \textbf{Indian paper} in  Bulletin of the Allahabad Mathematical Society paper from \textbf{2005} [25]. For GHW algebra setting and references see also an  \textbf{Indian paper} [20]. Finally - one may  contact   "`publikacje"' on http://ii.uwb.edu.pl/akk/publ1.htm .

\vspace{2mm}

\section{Snapshots on  spheres to comb and even parallelizable ones , seven "imaginary"  octonions'  units  in out of the Plato's Cave Reality quantum applications and 7Stones}

\vspace{2mm}

\noindent This is a section that should be written. Or rather it should be a book of books  telling the nowadays story of octonions' life in Mathematics, Gravity and Quantum theories  and invisibly playing a role in everyday's life of seven stones from the Cyber Space password http://www.7stones.com/.            

\noindent Here I am going to share with you some personal, hence partly  accidental expectations what should be there in this book of books. References quoted are by no means representative to the subject being defined on the way anyhow. Here come snapshots. Just like a running view through  the opened window of your running car. The weather  is all right. Snapshots.

\noindent \textbf{The Cauchy integral formulas and Cauchy - Riemann Equations}? For octonions?  - Yes. For example 
[26,27,...?...,13,14] and universal reference [5] with 94 references with  other references therein.

\vspace{2mm}
\noindent \textbf{Aspects of Quaternion and or Octonion Quantum Mechanics}? - Yes. For example [28-33] and universal reference source[5] with 94 references with  other references therein  and  [8] with week192. A contact with Paul C. Kainen's \textit{An Octonion Model for Physics} \textbf{2000} is recommended too [39]. The note  [30] is referring to  Octonionic Gauge Theory and the letter [31]  makes reference to seven dimensional octonionic Yang-Mills instantons and string solitons. A few months years old young and small note  [32]  brings us news about algorithms for least squares problem in quaternionic quantum theory.  Quaternionic Quantum Mechanics and Quantum Fields ? Yes and see why not?!  See [33]. It is the  608 pages book from year 1995.  Here we quote an abstract of its subjects' content:

\vspace{2mm}
\noindent "`It has been known since the 1930s that quantum mechanics can be formulated in quaternionic as well as complex Hilbert space. But systematic work on the quaternionic extension of standard quantum mechanics has scarcely begun. ... book signals a major conceptual advance and gives a detailed development and exposition of quaternionic quantum mechanics for the purpose of determining whether quaternionic Hilbert space is the appropriate arena for the long sought-after unification of the standard model forces with gravitation. Significant results from earlier literature, together with many new results obtained by the author, are integrated to give a coherent picture of the subject. The book also provides an introduction to the problem of formulating quantum field theories in quaternionic Hilbert space. The book concludes with a chapter devoted to discussions on where quaternionic quantum mechanics may fit into the physics of unification, experimental and measurement theory issues, and the many open questions that still challenge the field. ... ."'
\vspace{2mm}

\noindent   For those necessitating in further motive encouragement another quotation [from Girish Joshi's MATHEMATICAL STRUCTURES IN NATURE, Beyond Complex Numbers; just a glimpses:\\
http://www.ph.unimelb.edu.au/~ywong/poster/articles/joshi.html\\

\vspace{2mm}
\noindent "`IN NATURE, there is a deep connection between exceptional mathematical structures and the laws of micro- and macro-physics --- Quaternions and Octonions have played an important role in the recent development of pure and applied physics.

\noindent Quaternions were discovered by Hamilton in 1843 [19] and the Quaternions' main use in the 19th century consisted in expressing physical theories in "Quaternionic notation". An important work where this was done was Maxwell's treatise on electricity and magnetism. Toward the end of the century, the value of their use in electromagnetic theories led to a heated debate dubbed "The Great Quaternionic War".

\noindent In a 1936 paper, Birkhoff and von Neumann presented a propositional calculus for Quantum Mechanics and showed that a concrete realization leads to the general result that a Quantum Mechanical system may be represented as a vector space over the Real, Complex, and Quaternionic fields. Since then this area has remained active, aiming to extend Complex Quantum Mechanics (CQM) by generalizing the complex unit in CQM to Quaternions and to find observable effects of QQM. Jordan Algebras were proposed by Jordan, Neumann and Wigner in formulating non-associative Quantum Mechanics, where quantization is achieved through associators rather than commutators. This formulation allows mixing of space-time and internal symmetries. Another attractive feature of Jordan Algebra is that critical dimensions of 10 and 26 arise naturally, suggesting a connection to string theory.

\noindent Away from physics, Quaternions have recently been used for robotic control, computer graphics, vision theory, spacecraft orientation and geophysics. The space shuttle's flight software uses Quaternions in its guidance navigation and flight control computations.

\noindent  \textbf{"`Octonians"'.} "`In the early days of Quantum Mechanics, there were several attempts to introduce new algebraic structures in physics. In recent years, there has been remarkable activity in the field of physical applications of non-associative algebras: Octonian formulations of\textbf{ Yang}-Mills gauge theories and string theory, Octonian description of quarks and leptons, and Octonian supergravity."'

\vspace{2mm}
\noindent   At that hight of generality of a posteriori inclination Girish Joshi is perfectly right on the right side of the more and  more common convictions.

\vspace{2mm}

\noindent The Practical use of division algebras in gravity and cosmology Theories is not a new invention at all!  For example see the 1971 year paper [34].
For an  Online Article  posted on\textbf{ 2007-05-26} 20:16:11 abundant with ideas and links to highly  respective Mathemagics see [39].

\vspace{2mm}

\noindent \textbf{More on QM and octonions links  - another snapshots}

\noindent \textbf{Geoffrey Dixon}   Division Algebras: Octonions, Quaternions, Complex Numbers, and the Algebraic Design of Physics is at hand [9]. 
What is it about? \\ 
Abstract: "`The four real division algebras (reals, complexes, quaternions and octonions) are the most obvious signposts to a rich and intricate realm of select and beautiful mathematical structures. Using the new tool of adjoint division algebras, with respect to which the division algebras themselves appear in the role of spinor spaces, some of these structures are developed, including parallelizable 
spheres, exceptional Lie groups, and triality. In the case of triality the use of adjoint octonions greatly simplifies its investigation. Motivating this work, however, is a strong conviction that the design of our physical reality arises from this select mathematical realm. A compelling case for that conviction is presented, a derivation of the standard model of leptons and quarks."'

\vspace{2mm}

\noindent \textbf{Susumo Okubo}, Introduction to Octonion and Other Non-Associative Algebras in Physics. What is it about?
It covers topics ranging from algebras of observables in quantum mechanics and octonions to division algebra, triple-linear products and Yang–Baxter  equation. Here you may find the non-associative \textbf{gauge} theoretic reformulation of general relativity theory. [ You want to see present Tachion Kwa\'sniewski as a student talking with \textbf{Yang}? Chen Ning Yang? [\textbf{Yang}-Mills \textbf{gauge} theory] Contact my website http://ii.uwb.edu.pl/akk/index.html . In 1976 the present author was young and was invited by Yang to  SUNY - the  State Univesity of New York - for one semester. When I was ten -in 1957 Tsung Dao Lee and Chen Ning Yang won the Nobel Prize in physics.

\vspace{2mm}

\noindent \textbf{Murray Bremner } Quantum Octonions [29]. What is it about? Abstract: Bremner constructs a quantum deformation of the complex Cayley algebra  - a kind of  q-deformed version of the octonions. He  uses the representation  of $U_q (sl(2))$, the quantized enveloping algebra of the complex Lie algebra $sl(2)$.

\vspace{2mm}

\noindent  \textbf{Murat Gunaydin},  Hermann Nicolai,[31]; write on \textbf{ seven dimensional} octonionic Yang-Mills instanton.  What is it about?  Abstract: They construct an octonionic instanton solution to the\textbf{ seven }dimensional Yang-Mills theory based on the exceptional gauge group $G_2$ which is the automorphism group of this extreme  division algebra of octonions. This octonionic instanton has an extension to a solitonic solution of the low energy effective theory of the heterotic string that preserves two of the sixteen supersymmetries and hence corresponds to N = 1 space-time supersymmetry in (2+1) dimensions transverse to the \textbf{seven dimensions} where the \textbf{Yang}-Mills instanton is  defined. ...

\vspace{2mm}

\noindent \textbf{Seven Stones snapshot.} The Author:	Seven Stones Multimedia [7]. What is this about?  This is a series of tutorial modules which  apart from various vivid labs includes such topics as lessons about "`the \textbf{octonion algebra}, division algebras, lie algebras, lie groups, and spinors (Clifford algebras); introduction to special relativity; introduction to quantum mechanics, quantum mechanical uncertainty, nuclear decay and halflife lab, Bohr atom lab; statistics (introduction to linear regression); and astronomy (expanding universe, Doppler shift)."' 

\vspace{2mm}

\noindent \textbf{SO(8)}, \textbf{Triality}, \textbf{F4}, and \textbf{seven sphere  parallelizability ?}  . [One day one hade seen there in the sky over Poland  a flying exceptional structure  \textbf{F16}] with a new minister on board. This was not a Lie algebra, not at all. \\
Well, revenons a nos moutons! 
\noindent  Where you might be served the tutorial answer?  Naturally in [7] i.e in  \textbf{Seven Stones} Multimedia. Contact\\ http://www.7stones.com/Homepage/octotut14.html \\ i.e. the lesson 14 and others at wish. You might have a look also on  parallelizability  and triality writings [16,18] supposed to be concerned with physics theoretically.  
For more recent paper - see  [18]. Naturally  J. Baez is as always around and helpful [5,8].

\vspace{2mm}
\noindent The seven sphere  parallelizability  final history was born in $Cz\tilde{e}stochowa$ - the Holy Place of  Catholics of Poland and allover the Globe  being  so much praised by our  Santo subito John Paul II. [ Recall Castel Gandolfo meatings with Him - devoted also to Mathemagics.]\\
The final parallelizability   history was born there in $Cz\tilde{e}stochowa$ as far as the main contributor's birth is concerned.   

\noindent You know, the  classical problem was to determine which of the spheres $S^n$ are parallelizable. For $n = 1$ this is the circle in tha Gauss plane. For $n = 3$ the 3-sphere which is also the  group SU(2)and the only other parallelizable sphere is $S^7$. This was proved in 1958, by Michel André Kervaire, and  then by Raoul Bott and John Milnor, in independent work. note then  [see WIKI] he died on month ago:
Michel Andr\'e Kervaire  ($Cz\tilde{e}stochowa$, Poland, 26 April 1927 - Geneva, Switzerland,\textbf{ 19 November 2007})   
 
\noindent He worked for some time with \textbf{Milnor }[from SUNY - the  State Univesity of New York; Professor Yang was there too].  As \textbf{seven }is concerned  "`...John Milnor...-his most celebrated single result is his proof of the existence of \textbf{7-dimensional spheres }with nonstandard differential structure. Later with Michel Kervaire, he showed that the\textbf{ 7-sphere} has 15 differentiable structures (28 if you consider orientation). An n-sphere with nonstandard differential structure is called an exotic sphere, a term coined by Milnor."'\\
see http://mathworld.wolfram.com/Parallelizable.html.

\vspace{2mm}

\noindent Do not be afraid of nonassociativity of octonions. These are almost associative as forming  an example of the alternative algebra [13,14]. No surprise then that in [40] H. Albuquerque and S. Majid, observed that the "`octonions should have just the right properties as a substrate for the study of subtle properties of associativity in  tensor product algebras.   It is well-known that octaves form an {\it alternative} algebra; associativity holds when two of the three terms are equal; see, e.g., Schafer .  More; "`Albuquerque and Majid have proved that the octonions are associative up to a natural transformation"' - see also my invention in [13,14]. "`Their work uses ideas from quantum algebra and from Mac Lane's theory of coherence in categories.  The introduction of a thoroughly octaval viewpoint into the {\it topos} itself ought to have very interesting consequences for the enterprise of building a categorical model of continuum mechanics"'... [39].

\section{END}
\vspace{2mm}
\noindent \textbf{There is no end.}

\vspace{2mm}

\noindent \textbf{4.1.  Example  of no end.}

\noindent \textit{Added in proof} 1. 
After many whiles of good Google safe sale browsing  in the Internet Cosmos, after more readings and after receiving letters it is of 
outmost importance to indicate at least the following positions symptomatic for the subject  of our interest here.
This concerns supersymmetry  and the division algebras [41]  and integrable hierarchies [42].
We were also encouraged and inspired to look through interesting papers of Vladimir D. Dzhunushaliev thus
having been led to articles by Merab Gogberashvili [43,44].  The author of [43] uses the algebra of split octonions 
gaining an original octonionic version of Dirac Equations. The natural appearance of octonionic gradient gives us a hint to use
the octonionic analyticity as introduced by the author in [13,14]. The deeper gain seems to be achieved in the next paper [44]
on octonionic version of electrodynamics while noting that the non-existence of magnetic monopoles in classical electrodynamics
is related to an associativity limit  of Octonionic electrodynamics from [44].\\
Vladimir D. Dzhunushaliev works [45-49]  end up with "`Toy Models of a Nonassociative Quantum Mechanics"'[47] and 
supersymmetric extension of quantum mechanics via  introducing octonionios modulo quaternions! as "`hidden variables"'.
Both Merab Gogberashvili papers [43,44] as well as A. Dimakis and F. Muller-Hoissen [2] clue reference are referred  to by Vladimir D. Dzhunushaliev        
in his papers from our list of his references.\\
In order not to get an impression that OQM [ Octonionic Quantum Mechanics]  is a new subject let us at last but not least  recall as an example from the past century the important paper  [50] from 1996 by  Stefano De Leo  and  Khaled Abdel-Khalek  in which the authors  solve the hermiticity problem and define an appropriate momentum operator within OQM.
For more - including   also Octonionic Dirac Equation's matters see:
http://ptp.ipap.jp/cgi-bin/getarticle?magazine=PTP.....page=833-845 
There you find also substantial  references  to early papers by Professors  L. P. Horwitz and L. C. Biedenharn, [51] and Professor Jakub Rembielinski from Poland [52]. 
As a matter of fact the dissertation of  Professor Jakub Rembielinski was my first contact with
OQM.  Of course apart from fundamental [no way out?] paper from  1934 by  Jordan , von Neumann and Wigner [53] which I had gone through as a student http://ii.uwb.edu.pl/akk/.

\noindent \textbf{4.2.  Example  of no end.}

\noindent The recent references  from People's Republic of China on quaternion's applications that follow  those from the \textbf{4.1.  Example  of no end.} are the effect of e-mail correspondence.

\noindent \noindent \textbf{4.3.  Example  of no end.}
\noindent The Reader is welcomed to contribute. The Example is opened.

\vspace{5mm}

\noindent \textbf{Acknowledgements}
The author  thanks Maciej Dziemia\'nczuk,  the Student of the Institute of Computer Science at Bialystok University for his kind interest to this work and also for being helpful with LaTeX matters. Thanks and Glory Unto Him. I express my gratitude to Tevian Dray and Joseph Malkevitch  for permission to use their fine Fano plane octonions pictures as well as for Broome Bridge reproduction of Tevian Dray.  This essay is \textbf{dedicated} to my sons :
the Master of computer science and practice - Jakub the First , the gifted Mathematician Bartosz-Kosma\\
http://math.uwb.edu.pl/~zaf/kwasniewski/index.html
and Philosopher of Art and Eternal Life Kosma-Damian.
And the last but not least thanks to my wife for being so withstanding person during the past and forthcoming life keeping  calm in the midst of the storm of  mathemagics and philosophy of her sons and husband.

\end{document}